\documentclass[10pt]{amsart}

\usepackage{amsmath,amssymb,amsfonts,epsfig}

\newcommand{\rb}{\raisebox}
\newcommand{\ig}{\includegraphics}
\newcommand\risS[6]{\rb{#1pt}[#5pt][#6pt]{\begin{picture}(#4,15)(0,0)
  \put(0,0){\ig[width=#4pt]{#2.eps}} #3
     \end{picture}}}
\def\a{\alpha}
\def\b{\beta}
\def\d{\delta}
\def\cS{\mathcal S}
\newtheorem{thm}{Theorem}[section]
\newtheorem{defn}[thm]{Definition}

\def\kb#1{[ #1 ]}

\def\R{{\mathbb R}}
\begin{document}

\title[Bollob\'as-Riordan and relative Tutte polynomials]
  {Bollob\'as-Riordan and relative Tutte polynomials}
\author{CLARK BUTLER and SERGEI CHMUTOV}
\date{}

\address{Department of Mathematics, The Ohio State University,
231 West 18th Avenue, Columbus, OH 43210\linebreak 
{\tt butler.552@buckeyemail.osu.edu}\linebreak 
{\tt chmutov@math.ohio-state.edu}}


\begin{abstract}
We establish a relation between the Bollob\'as-Riordan polynomial of a ribbon graph with the relative Tutte polynomial of a plane graph obtained from the ribbon graph using its projection to the plane in a nontrivial way. Also we give a duality formula for the relative Tutte polynomial of dual plane graphs and an expression of the Kauffman bracket of a virtual link as a specialization of the relative Tutte polynomial.
\end{abstract}

\maketitle

\section*{Introduction} \label{s:intro}
Graphs on surfaces can be studied in terms of plane graphs via a projection preserving the rotation system. These plane graphs are constructed in a nontrivial way in order to preserve the topological information of the graph. These constructed plane graphs usually have some extra (distinguished) edges and extra vertices. They are {\em relative} plane graphs.

\bigskip
\begin{defn}
A {\it relative plane graph} is a plane graph $G$ with a distinguished subset 
$H \subseteq E(G)$ of edges. The edges $H$ are called the 0-{\it edges} of $G$.
\end{defn}

For example, the presence of artificial virtual crossings, as in the picture of $K_5$ below, may be captured by 0-edges and some additional vertices. 
$$\risS{-40}{K5-vir-cr}{\put(-30,-25){\parbox{2.5in}{
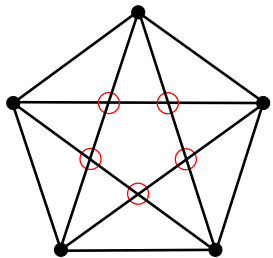
          \tt Virtual crossings in plane pictures of non-planar graphs}}
                       }{100}{60}{80}\hspace{1in}
\risS{-40}{K5-rel-plan}{\put(-10,-25){\parbox{2.3in}{
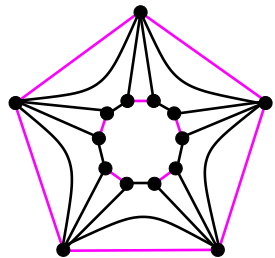
          \tt Corresponding relative plane graph; 0-edges are dark}}
                       }{100}{60}{0}
$$

The motivation of our work comes from the knot theory. The classical Thistlethwaite theorem \cite{Th} relates the Jones polynomial of an alternating link to the Tutte polynomial of plane graph obtained form a checkerboard coloring of the regions of the link diagram. This theorem has two different kinds of generalizations to virtual links. 
One \cite{Ch, CP, ChVo, DFKLS, Mof} involves graphs on surfaces and a topological version of the Tutte polynomial due  to B.~Bollob\'as and O.~Riordan \cite {BR3}. Another generalization is based on a relative version of the Tutte polynomial found by Y.~Diao and G.~Hetyei \cite{DH}. 
In this paper we establish a direct relation between the Bollob\'as--Riordan and relative Tutte polynomials that explains how these two generalizations are connected.

In Section \ref{s:constr} we explain the construction of a relative plane graph from a ribbon graph as well as how to recover a ribbon graph from a relative plane graph.
Our main theorem is formulated in Section \ref{s:th} and proved in Section
\ref{s:proof}. In Section \ref{s:dual} we describe the relation between
relative Tutte polynomials of dual plane graphs generalizing the classical relation $T_{G}(x,y)=T_{G^*}(y,x)$. In Section \ref{s:Kabr} we obtain the Kauffman bracket of a virtual link in terms of the relative Tutte polynomial, improving the theorem of \cite{DH}.

This work has been done as a part of the Summer 2010 undergraduate research working group
\begin{center}\verb#http://www.math.ohio-state.edu/~chmutov/wor-gr-su10/wor-gr.htm#
\end{center}
``Knots and Graphs" 
at the Ohio State University. We are grateful to all participants of the group for valuable discussions.

\section{Ribbon graphs and relative plane graphs} \label{s:constr}

By a {\it ribbon graph} we mean an abstract (not necessarily orientable) surface with boundary decomposed into topological discs of two types, {\it vertex-discs} and {\it edge-ribbons}, satisfying some natural conditions; see \cite{BR3,Ch} for precise definition and 
\cite{GT,LZ,MT} for the general notions and terminology of topological graph theory. Ribbon graphs are considered up to homeomorphisms of the underlying surfaces preserving the decomposition. A ribbon graph can be regarded as a regular neighborhood of a graph cellularly embedded into a surface. Thus the language of ribbon graph is essentially the same as for cellularly embedded graphs.

\subsection{From ribbon graphs to relative plane graphs.}\label{ss:rg2rpg}\ 

Let $R$ be a ribbon graph. 
Consider a projection $\pi: R\to \R^2$ which is 1-to-1 except the points of singularities. We restrict these singularities to two types. The first occurs when a ribbon twists over itself; in this case a whole line interval on the ribbon is projected to a single point. The second type occurs when the images of two edge ribbons cross. In this case, the projection is 2-to-1 over the disc of the intersection. 
$$\risS{-20}{singular}{\put(60,50){$R$}\put(30,30){$\pi$}
         }{60}{45}{45}\hspace{3cm}
\risS{-30}{singular2}{\put(45,60){$R$}\put(25,40){$\pi$}
        \put(-150,-8){\tt Possible singularities of the projection}
         }{40}{0}{0}\hspace{2cm}
$$
Consequently, the restriction of the projection to the boundary of the ribbon graph is am immersion with only double points as singularities.

The construction of the relative plane graph $G$ consists of four steps.
\begin{enumerate}
\item
On each edge of $R$ we choose a portion of the ribbon on which the projection is 1-to-1. We call it and its image on $\R^2$ a {\it regular edge}. The regular edges will be the non-zero edges of the relative plane graph $G$.
\item
Extend the vertex discs of $R$ through to the regular edge of each ribbon.
\item
Each of these extended vertices is segmented by the regular edges and the singularities of the projection. These segments become the vertices of $G$. 
\item
The {\it 0-edges} of $G$ correspond to the double points of the restriction of the projection to the boundary of $R$. They connect the vertices of $G$ which correspond to the extended regions sharing the same double point in a checkerboard manner.
\end{enumerate}

{\bf Examples}.\\
{\bf 1)} The picture of the graph $K_5$ from the introduction gives such a projection. The next figure shows the steps in construction of the corresponding graph $G$
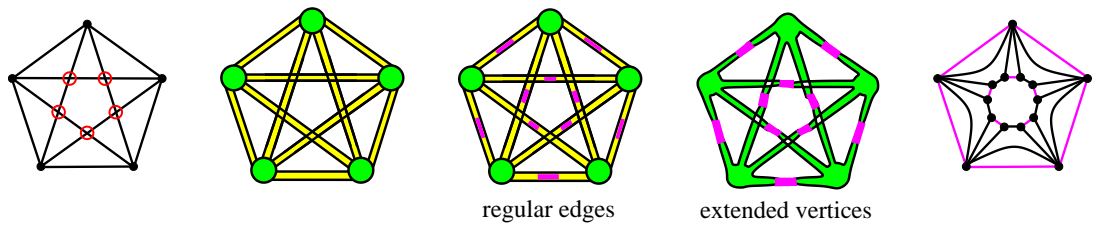
$$\risS{-20}{K5-vir-cr-sm}{}{60}{50}{40}\qquad
\risS{-25}{Ga-K5}{}{70}{0}{0}\qquad
\risS{-40}{Ga-K5-re}{}{70}{0}{0}\qquad
\risS{-38}{Ga-K5-ev}{}{70}{0}{0}\qquad
\risS{-20}{K5-rel-plan-sm}{}{60}{0}{0}
$$

Here are two more examples showing the dependence of $G$ on the choice of regular edges.

{\bf 2)}
$$\hspace{-50pt}\risS{-100}{rg1a}{}{350}{50}{120}$$

{\bf 3)}
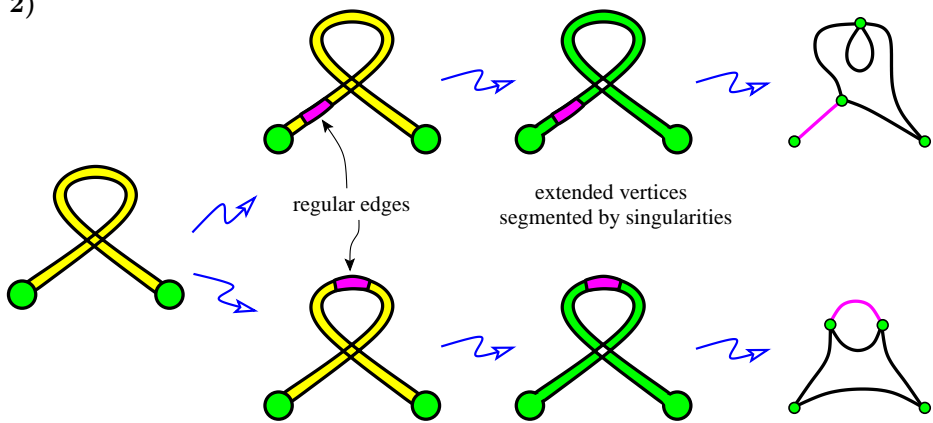
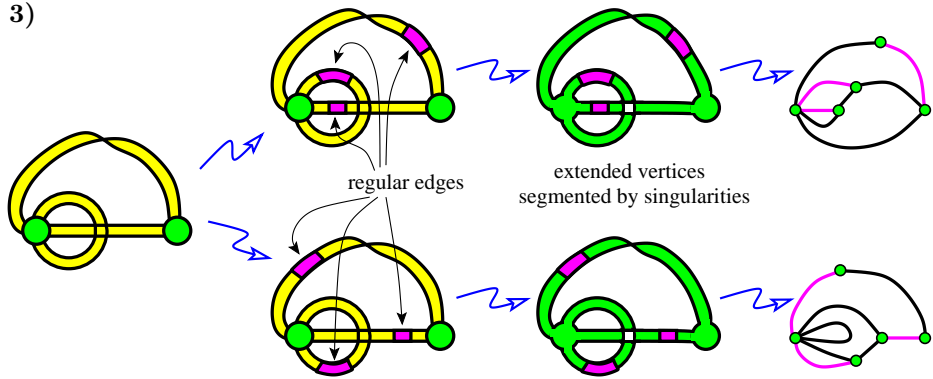
$$\hspace{-50pt}\risS{-100}{rg2a}{}{350}{30}{110}$$
It is a remarkable consequence of the main theorem that the specialization of the relative Tutte polynomial does not depend on all these choices. It  is not difficult to describe a sequence of moves on relative plane graphs relating the graphs with different choices of the regular edges. It would be interesting to find such moves for different choices of the projection $\pi$ and, more generally, the moves preserving the relative Tutte polynomial.

\subsection{From relative plane graphs to ribbon graphs.}\label{ss:rpg2rg}\ 

Conversely, from a relative plane graph $G$ we may construct a ribbon graph $R$ in the following way.
\begin{enumerate}
\item
Delete all regular edges of $G$ and for the obtained plane graph $H$ construct its medial graph.
\item
Consider the medial graph as an immersion of a union of circles with double points on its vertices, midpoints of 0-edges of $G$. Each such circle gives a {\it vertex} of $R$. 
\item 
Each regular edge of $G\setminus H$ makes a trace on the circles which we denote by two arrows according to a counterclockwise orientation of the plane.
\item
Pull apart and untwist the circles. Here we may reverse the orientations of certain arrows.
\item
The {\it edges} of $R$ are formed by attaching ribbons to the vertex discs of $R$ according to the arrows. Thus we get the arrow presentation of the ribbon graph $R$ in sense of \cite{Ch}.
\end{enumerate}

\bigskip
{\bf Examples}.\\
{\bf 1)} 
$$\risS{-100}{pg-rg1}{}{350}{10}{100}$$
{\bf 2)}
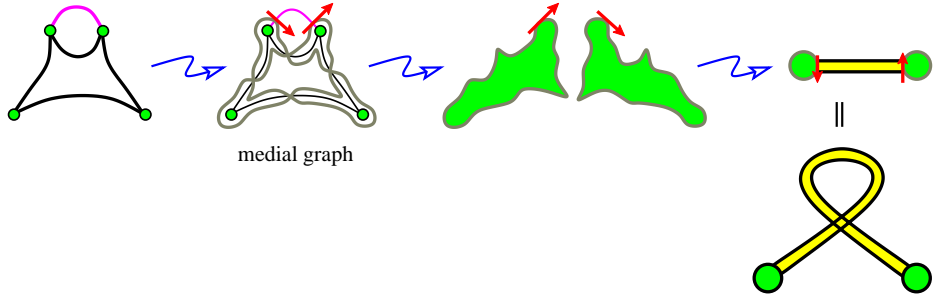
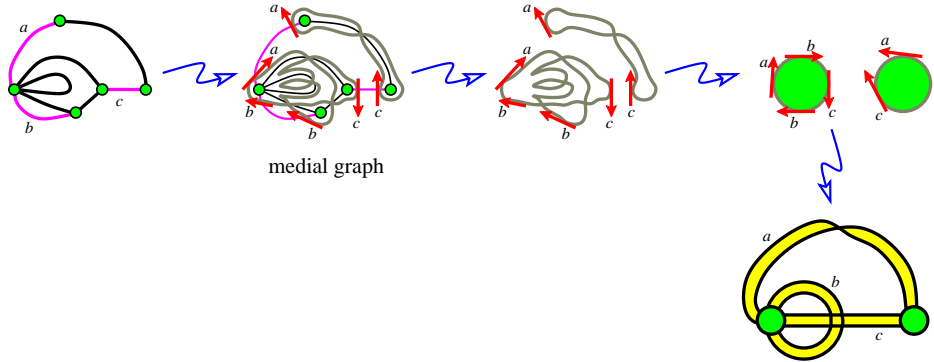
$$\risS{-100}{pg-rg2}{}{350}{40}{110}$$
Observe that this construction depends on the planar embedding of $G$, as different embeddings can change the cyclic order of edges on a vertex of $R$. However, if we fix a projection of $R$, the regular edge portions of ribbons of $R$, and the planar embedding of $G$ arising from this projection, then this construction is inverse to the construction of a relative plane graph from a ribbon graph. 

\subsection{The Bollob\'as--Riordan polynomial of ribbon graphs.}\ 

The Bollob\'as--Riordan polynomial, originally defined in \cite{BR3}, was generalized to a multivariable polynomial of weighted ribbon graphs in
\cite{Mof, VT}. We will use a sightly more general doubly weighted 
Bollob\'as--Riordan polynomial of a ribbon graph $R$ with weights $(x_e,y_e)$ of an edge $e\in R$.

\bigskip
\begin{defn}
$$B_{R}(X,Y,Z):=\sum_{F\subseteq R} (\prod_{e\in F}x_e) 
(\prod_{e\in R \setminus F}y_e) X^{k(F)-k(R)}Y^{n(F)}Z^{k(F)-bc(F)+n(F)}\ ,
$$
where the sum runs over all spanning subraphs $F$, $k(F)$ is the number of connected components of $F$, $n(F)=|E(F)|-v(F)+k(F)$ is the {\it nullity} of $F$, and $bc(F)$ is the number of boundary components of $F$. 
\end{defn}

\subsection{The relative Tutte polynomial.}\ 

\bigskip
\begin{defn}
 Let $G$ be a relative plane graph with the distinguished set of 0-edges $H$.
$$T_{G,H}(G):=\sum_{F \subseteq G\setminus H} (\prod_{e \in F}x_e) 
(\prod_{e \in \overline{F}}y_e) X^{k(F\cup H)-k(G)}Y^{n(F)}\psi (H_F)\ ,
$$
where $\overline{F}=G \setminus (F\cup H)$, $\psi$ is a block-invariant function on graphs, and $H_F$ is the plane graph obtained from $F \cup H$ by contracting all edges of $F$. Our choice of $\psi$ is 
$$\psi (H_F):=d^{\delta(H_F)-k(H_F)}w^{v(H_F)-k(H_F)}\ ,
$$
$\delta(H_F)$ is the number of circles that immerse to the medial graph of $H_F$. 
\end{defn}

\bigskip
{\bf Remarks.}\ 

{\bf 1.} The relative Tutte polynomial was introduced by Y.~Diao and G.~Hetyei in \cite{DH}, who use the notion of {\it activities} to produce the most general form of it. The all subset formula we use was discovered by a group of undergraduate students (M.~Carnovale, Y.~Dong, J.~Jeffries) at the  OSU summer program ``Knots and Graphs" in 2009.
However, the similar expressions may be traced back to L.~Traldi \cite{Tr}
for non relative case, and to S.~Chaiken \cite{Cha} for relative case of matroids. 

{\bf 2.} The function $\psi$ in \cite{DH} can be obtained from ours by substitution $w=1$. 

{\bf 3.} Another difference with \cite{DH} is that we are using a doubly weighted version of the relative Tutte polynomial with weights $(x_e,y_e)$ of an edge $e\in G\setminus H$.

{\bf 4.} In the process of constructing the graph $H_F$ by contracting the edges of $F$ in $F \cup H$, we may come to a situation when we have to contract a loop. Then under contraction of a loop we actually mean its deletion. Since $G$ and $F \cup H$ are plane graphs, then the graph $H_F$ is also embedded into the plane.

\section{Main Theorem} \label{s:th}

\begin{thm}
{\it Suppose $R$ is a ribbon graph, and $G$ is a relative plane graph of a projection of $R$. Or, equivalently, $G$ is a relative plane graph and $R$ is the ribbon graph arising from $G$. 

Then under the substitution $w=\sqrt{\frac{X}{Y}}, d=\sqrt{XY}$,
$$X^\a Y^\b T_{G,H}(X,Y)=B_R ( X,Y,\frac{1}{\sqrt{XY}} )\ ,
$$
where $\a:=k(G)-k(R)-\b$ and $\b:=-\frac{1}{2}(v(R)-v(G))$.}
\end{thm}

\bigskip
{\bf Remark.} The construction of $G$ from $R$ and backward can be generalized to a wider class of projections $\pi$. We can require only the restriction of $\pi$ to the boundary of $R$ to be an immersion with only  ordinary double points as singularities. The theorem holds in this topologically more general situation. However, from the point of view of graph theory it is more natural to restrict ourselves to the class of projections which we use, as the more general projections are not necessarily invertible. 

\section{Proof} \label{s:proof}

Our constructions of $G$ from $R$ and backward from Sections \ref{ss:rg2rpg} and \ref{ss:rpg2rg} give a bijection between (regular) non-zero
edges of $G$ and the edge-ribbons of $R$. We denote the corresponding edges by the same letter $e$ for both $e\in G\setminus H$ and for $e\in R$ since this will not lead to confusion. Moreover, in the theorem we assume that this bijection respects the weights of the doubly weighted polynomials.
The bijection can be naturally extended to the bijection between spanning subgraphs $F\subseteq G\setminus H$ and 
$F'\subseteq R$ so that the weights of $F$ and $F'$ are equal to each other:
$$(\prod_{e \in F}x_e) (\prod_{e \in \overline{F}}y_e) =
(\prod_{e\in F'}x_e) (\prod_{e\in R \setminus F'}y_e)
$$

Thus the theorem can be checked only on monomials in $X$ and $Y$ corresponding to $F\subseteq G\setminus H$ and $F'\subseteq R$. In other words, we have to prove that
\begin{multline}\label{eq:mon}
X^{k(G)-k(R)--\frac{1}{2}(v(R)-v(G))} Y^{-\frac{1}{2}(v(R)-v(G))} 
X^{k(F\cup H)-k(G)}Y^{n(F)}
d^{\delta(H_F)-k(H_F)}w^{v(H_F)-k(H_F)} \\ 
=X^{k(F')-k(R)}Y^{n(F')}Z^{k(F')-bc(F')+n(F')}
\end{multline}
for $d=\sqrt{XY}$, $w=\sqrt{\frac{X}{Y}}$, and $Z=\frac{1}{\sqrt{XY}}$.

We need the following combinatorial equalities:
\begin{enumerate}
\item[(2)] $|E(F)|=|E(F')|$
\item[(3)] $k(H_F)=k(F \cup H)$
\item[(4)] $bc(F')=n(F)+\delta (H_F)$
\item[(5)] $v(H_F)=k(F)$
\end{enumerate} 

(2) is clear from the subgraph correspondence. 
Since contracting edges of a graph cannot disconnect it or connect disconnected components, (3) is immediate.

\subsection{Proof of (4).}\ 

In the plane subgraph $F \cup H$ of $G$, trace the edges two times each in the following manner: the edges of $F$ are traced by parallel lines and the edges of $H$ are traced by crossed lines as in the medial graph:
$$F\cup H = \risS{-28}{FcupH}{\put(-60,45){$F=\{a,b\}$}
                               \put(195,28){$F'\!\!=$}}{300}{25}{30}
$$
The number of immersed circles of the tracing is precisely $bc(F')$, by our construction.  

Now consider contracting edges of $F$. The contraction of a non-loop does not change the number of curves of the tracing. 
$$\risS{0}{contr-nl}{}{250}{40}{0}
$$

However, the contraction of a loop, i.e. deletion of the loop, fuses two disjoint curves together, one from the outside of the loop and one from the inside of the loop. So it reduces the number of curves by 1. 
$$\risS{0}{contr-l}{}{350}{50}{0}
$$

The result of contracting all the edges of $F$ is the graph $H_F$, for which the number of curves will be $\d(H_F)$. Since the number of loops contracted is $n(F)$, we have
$$bc(F')=n(F)+\d(H_F)\ .
$$

\subsection{Proof of (5).}\ 

Let us consider $F \cup H$ as a spanning subgraph of $G$ and let us remove the edges of
$H$ from it for a moment. Then we get the spanning subgraph $F$. Its edges are supposed to be contracted, so each connected component of $F$ gives a vertex of the resulting graph. Now restoring the edges of $H$ does not change the number of vertices of graph obtained by contracting by $F$. Thus $v(H_F)=k(F)$.
%
%

\subsection{Proof of the theorem.}\ 

We deal with the exponents of $X$ and $Y$ separately. The exponent of $X$ in the left hand side of equation \eqref{eq:mon} is 
$$\frac{1}{2} (v(R)-v(G))+k(G)-k(R)+k(F \cup H) - k(G) + \frac{1}{2}(\d(H_F)-k(H_F)+v(H_F)-k(H_F))
$$

Substituting the equalities above and making appropriate cancellations, 
\[\begin{split}
& =\frac{1}{2}(v(R)-v(G))-k(R) + \frac{1}{2}(\d(H_F)+k(F))\\
& =\frac{1}{2}(v(R)-v(G))-k(R) + \frac{1}{2}(bc(F')-n(F)+k(F))\\
& =\frac{1}{2}(v(R)-v(G))-k(R) + \frac{1}{2}(bc(F')-|E(F')|+v(G))\\
& = -k(R) + \frac{1}{2}(bc(F') + v(R) - |E(F')|)\\
& = k(F') - k(R) + \frac{1}{2}(bc(F')-n(F')+k(F'))\\
\end{split}
\]
which is the exponent of $X$ in $B_R$.

\bigskip
For $Y$, the exponent in the left hand side of equation \eqref{eq:mon} is
$$-\frac{1}{2}(v(R)-v(G)) + n(F) + \frac{1}{2}(\d(H_F)-k(H_F)-v(H_F)+k(H_F))
$$

This is equival to,
\[\begin{split}
& = |E(F)|-v(G)+k(F) + \frac{1}{2}(bc(F')-n(F)-v(H_{F})-v(R)+v(G))\\
& = \frac{1}{2}(bc(F')-|E(F)|+2v(G)-v(R)-2k(F))+|E(F)|-v(G)+k(F)\\
& = \frac{1}{2}(|E(F')|-v(R)+bc(F'))\\
& = n(F') + \frac{1}{2}(bc(F')-n(F')-k(F'))
\end{split}
\]
which is the exponent of $Y$ in $B_R$.

\section{Dual relative plane graphs} \label{s:dual}

Let $G$ be a relative plane graph. For an embedding of $G$ in the plane, the {\it dual} of $G$, denoted $G^*$ is formed by taking the dual of $G$ as a plane graph, and labeling the edges of $G^*$ which intersect 0-edges of $G$ as the 0-edges of $G^*$. Note that $(G^*)^*=G$, as with usual planar duality.

\begin{thm}
Under the substitution $w=\sqrt{\frac{X}{Y}}$, $d=\sqrt{XY}$, we have
\[
X^{a(G,H)}Y^{b(G)}T_{G,H}(X,Y)=Y^{a(G^*,H^*)}X^{b(G^*)}T_{G^*,H^*}(Y,X)
\]
with the correspondence on the edge weights being $x_{e}$=$y_{e^*}$, $y_{e}$=$x_{e^*}$, where $e^*$ is the edge of $G^*$ that intersects $e$, and\qquad
$a(G,H) = (|E(G \setminus H)|-v(G))/2+k(G)\ ,\qquad
b(G) = v(G)/2\ .$
\end{thm}

{\bf Proof of the Theorem.} 
The equality is on monomials of $T_{G,H}$, $T_{G^*,H^*}$ in the edge weights variables $(x_e,y_e)$ which establishes the correspondence between spanning subgraphs $F$ of $G\setminus H$ and $F^*$ of $G^*\setminus H^*$. Namely, $F^*$ consists of those regular edges of $G^*$ which do not intersect the regular edges of $F$.
 
We prove the equality on monomials for the exponent of $X$. Equality for $Y$ then follows from duality. The exponent of X on the left is
\[
\begin{split}
& \frac{1}{2}(|E(G \setminus H)|-v(G))+k(G)+k(F \cup H)-k(G)+ \frac{1}{2}(\d(H_{F})-k(H_{F})+v(H_{F})-k(H_{F})) \\
&\hspace{100pt} = \frac{1}{2}(|E(G \setminus H)|-v(G)+bc(F_R)-n(F)+k(F))\\
&\hspace{100pt}  = \frac{1}{2}(|E(G \setminus H)|+bc(F_R)-|E(F)|)\\
&\hspace{100pt}  = \frac{1}{2}(|E(\overline{F})|+bc(F_{R}))
\end{split}
\]
where $F_R$ is the ribbon graph constructed from the relative plane graph $F \cup H$ in the manner of Section \ref{ss:rpg2rg}.

On the right, let $F^*$ denote the subgraph of $G^*$ corresponding to $F$. Then the exponent of $X$ is\vspace{-8pt}
\[
\begin{split}
& n(F^*)+\frac{1}{2}(\d(H_{F^*})-k(H_{F^*})-v(H_{F^*})+k(H_{F^*})+v(G^*)) \\
& = n(F^*)+\frac{1}{2}(bc(F_R^*)-n(F^*)-k(F^*)+v(G^*)) \\
& = \frac{1}{2}(bc(F_R^*)+|E(F^*)|)
\end{split}
\]
Now, $|E(F^*)|=|E(\overline{F})|$ by the subgraph correspondence. The equality $bc(F_R)$=$bc(F_R^*)$ follows from the fact that the ribbon graphs $F_R$ and $F_R^*$ have the same boundary. Also it can be seen from the following figures:

$$\risS{0}{dual-1}{\put(35,31){$e$}\put(15,-5){$G\ni e\not\in F$}
              \put(146,26){$e^*$}\put(108,-12){$G^*\ni e^*\in F^*$}
                  }{180}{40}{20} \hspace{2cm}
\risS{0}{dual-1a}{\put(35,35){$e$}\put(15,-5){$G\ni e\in F$}
              \put(142,26){$e^*$}\put(108,-12){$G^*\ni e^*\not\in F^*$}
                  }{180}{40}{20}
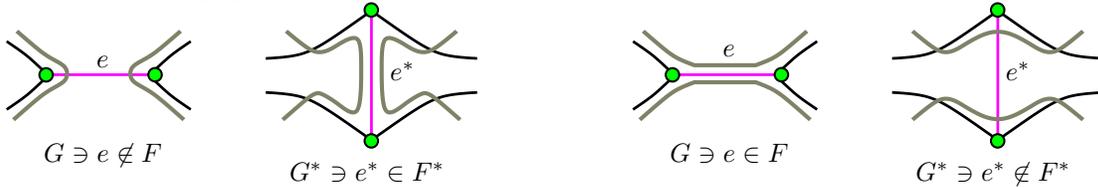
$$

\section{Kauffman bracket of virtual links} \label{s:Kabr}

In this section we generalize the result of \cite{DH} which extends the 
Thistlethwaite theorem to virtual links. 
Virtual links are represented by diagrams similar to ordinary knot diagrams, except some crossings are designated as {\it virtual}. Here are some examples of virtual knots.
$$\risS{-18}{ex}{}{65}{20}{20}\hspace{3cm}
  \risS{-18}{v31}{}{40}{0}{20}\hspace{3cm}
  \risS{-18}{v41}{}{55}{0}{20}\label{ex:vir-kn}
$$

Virtual link diagrams are considered up to plane isotopy, the {\it classical}
Reidemeister moves:
$$\risS{-10}{RI}{}{65}{20}{12}\qquad\qquad
  \risS{-10}{RII}{}{65}{0}{0}\qquad\qquad
  \risS{-10}{RIII}{}{65}{0}{0}\quad ,
$$
and the {\it virtual} Reidemeister moves:
$$\risS{-10}{RI-v}{}{65}{17}{15}\qquad\quad
  \risS{-10}{RII-v}{}{65}{0}{0}\qquad\quad
  \risS{-10}{RIII-v}{}{65}{0}{0}\qquad\quad
  \risS{-10}{RIV-v}{}{65}{0}{0}\quad .
$$

The Kauffman bracket for virtual links is defined in the same way as for classical links. Let $L$ be a virtual link diagram.
Consider two ways of resolving a classical crossing.
The {\it $A$-splitting},\ 
$\risS{-4}{cr}{}{15}{15}{8}\ \leadsto\ \risS{-4}{Asp}{}{15}{0}{0},$
is obtained by joining the two vertical angles swept out by the overcrossing arc when
it is rotated counterclockwise toward the undercrossing arc.
Similarly, the {\it $B$-splitting},\ 
$\risS{-4}{cr}{}{15}{15}{8}\ \leadsto\ \risS{-4}{Bsp}{}{15}{0}{0},$
is obtained by joining the other two vertical angles. A {\it state} $s$ of
a link diagram $L$
is a choice of either an $A$ or $B$-splitting at each classical crossing.
Denote by $\cS(L)$ the set of states of $L$.
A diagram $L$ with $n$ crossings has $|\cS(L)| = 2^n$
different states.

Denote by $\a(s)$ and $\b(s)$ the numbers of $A$-splittings and $B$-splittings
in a state $s$, respectively, and by $\d(s)$ the number of
components of the curve obtained from the link
diagram $L$ by 
splitting according to the state $s \in \cS(L)$. Note that virtual crossings do not connect components.

\begin{defn}\label{def:kb}
The \emph{Kauffman bracket} of a diagram $L$ is a polynomial in three variables
$A$, $B$, $d$ defined by the formula
$$ \kb{L} (A,B,d)\ :=\ \sum_{s \in \cS(L)} \,
A^{\a(s)} \, B^{\b(s)} \, d^{\,\d(s)-1}\,.
$$
\end{defn}

Note that $\kb{L}$ is \emph{not} a topological invariant of the link; it depends on the link diagram and changes with Reidemeister moves. However, it determines the \emph{Jones polynomial}
$J_L(t)$ by a simple substitution:\vspace{-5pt}
$$A=t^{-1/4},\qquad B=t^{1/4},\qquad d=-t^{1/2}-t^{-1/2}\ ;$$
$$J_L(t)\, := (-1)^{w(L)} t^{3w(L)/4} \kb{L} (t^{-1/4}, t^{1/4}, -t^{1/2}-t^{-1/2})\ .
$$


In 1987 Thistlethwaite \cite{Th} (see also \cite{Ka1}) proved that up to a sign and a power of $t$ the Jones polynomial 
$V_L(t)$ of an alternating link $L$ is equal to 
the Tutte polynomial $T_{G_L}(-t,-t^{-1})$
of the Tait graph $G_L$ obtained from a checkerboard coloring of the
regions of a link diagram. 
$$\risS{-7}{kd-gr}{\put(0,40){\mbox{$L$}}  \put(190,40){\mbox{$G_L$}}
                 }{200}{45}{10} 
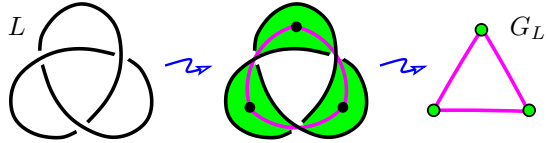
\label{p:thist}
$$
L.~Kauffman \cite{Ka2} generalized the theorem to arbitrary (classical) links using signed graphs. To virtual links this theorem was extended in \cite{Ch, CP, ChVo} using ribbon graphs. Another extension, using the relative Tutte polynomial, is due Y.~Diao and G.~Hetyei \cite{DH}. In their construction the relative plane graph is the Tait graph of a virtual link diagram whose 0-edges correspond to virtual crossings. They expressed $\kb{L}(A,A^{-1},-A^2-A^{-2})$ as a specialization of the relative Tutte polynomial. 
We think that the whole Kauffman bracket $\kb{L}(A,B,d)$, although not a link invariant, is of some interest as a pure combinatorial invariant of link diagrams. It turns out that it also can be expressed as a specialization of the relative Tutte polynomial.

\bigskip
Following \cite{DH}, we assign signs to the edges of the Tait graph $G$
depending on whether the edge connects  $A$- or $B$-splitting regions: \qquad
$\risS{-30}{cr_pl}{}{35}{0}{30}\hspace{2cm}
  \risS{-30}{cr_mi}{}{35}{0}{0}$

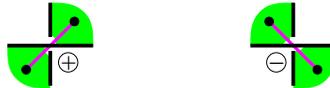
\begin{thm}
{\it Let L be a virtual link diagram, and G the relative plane Tait graph of L. Then, under the substitution 
$$X=\frac{Bd}{A},\qquad Y=\frac{Ad}{B},\qquad w=\frac{B}{A},\qquad x_+=y_+=1,\qquad  x_-=\sqrt{\frac{X}{Y}}\!=\!\frac{B}{A},\qquad y_-=\sqrt{\frac{Y}{X}}\!=\!\frac{A}{B}
$$
we have,
$$[L](A,B,d)=A^{v(G)-k(G)}B^{|E(G \backslash H)|-v(G)+k(G)}d^{k(G)-1}T_{G, H}\ .
$$}
\end{thm}

{\bf Proof.} 

The equality is on monomials, with the correspondence between subgraphs $F$ and states $S$ being the natural one: 
$$\risS{-30}{F2S}{\put(10,43){$e$}\put(78,68){$A$}\put(78,6){$B$}
        \put(30,54){\rotatebox{15}{$e\in F$}}
        \put(28,23){\rotatebox{-25}{$e\not\in F$}}
        \put(177,43){$e$}\put(236,68){$A$}\put(236,6){$B$}
        \put(186,53){\rotatebox{15}{$e\not\in F$}}
        \put(188,23){\rotatebox{-25}{$e\in F$}}
              }{250}{40}{30}
$$
Let $|E_-(F)|$ (resp. $|E_+(F)|$) be the number of negative 
(resp. positive) edges in the graph $F$. The power of $B$ on the right is
\[\begin{split}
& |E(G \setminus H)|-v(G)+k(G)+|E_-(F)|-|E_-(\overline{F})|+k(F \cup H)-k(G)                                    -n(F)+ v(H_F) - k(H_F) \\
& = |E_-(F)|-|E_-(\overline{F})|+|E(G \setminus H)|-|E(F)|\\
& = |E_-(F)|-|E_-(\overline{F})|+|E(\overline{F})|\ =\ 
     |E_-(F)|+|E_+(\overline{F})|\ =\  \b(S)\ ,
\end{split}
\]
as it can be easily seen from the picture above.
The proof of equality on the exponent of $A$ is similar. 
For $d$, the exponent on the right is
$$k(G)-1+k(F \cup H)-k(G)+n(F)+\d(H_F)-k(H_F) = n(F)+\d(H_F)-1
 = bc(F_R)-1 = \d(S)-1\ .
$$


\bigskip
\begin{thebibliography}{ABCDE}

\bibitem[B]{B} B.~Bollob\'as, {\it Modern graph theory},
   Graduate Texts in Mathematics~{\bf 184}, Springer, New York, 1998.

\bibitem[BR]{BR3} B.~Bollob\'as and O.~Riordan, {\it A polynomial of graphs
   on surfaces}, Math.~Ann.~{\bf 323} (2002) 81--96.

\bibitem[Cha]{Cha} S.~Chaiken,  {\it The Tutte polynomial of a ported 
   matroid}, 
   Journal of Combinatorial Theory, Ser. B, {\bf 46} (1989) 96--117.

\bibitem[Ch]{Ch} S.~Chmutov, {\it Generalized duality for graphs on 
   surfaces and the signed Bollob\'as-Riordan polynomial}, 
   Journal of Combinatorial Theory, Ser. B, {\bf 99}(3) (2009) 617--638;
    preprint \verb#arXiv:math.CO/0711.3490#, 

\bibitem[ChPa]{CP} S.~Chmutov, I.~Pak,  {\it The Kauffman bracket of 
    virtual links and the Bollob\'as-Riordan polynomial}
    Moscow Mathematical Journal~{\bf 7}(3) (2007) 409--418; 
    preprint \verb#arXiv:math.GT/0609012#, 
   
\bibitem[ChVo]{ChVo} S.~Chmutov, J.~Voltz,  {\it Thistlethwaite's theorem 
   for virtual links}, 
   Journal of Knot Theory and Its Ramifications, {\bf 17}(10) (2008)  
   1189–-1198; preprint \verb#arXiv:math.GT/0704.1310#.

\bibitem[DFKLS]{DFKLS} O.~Dasbach, D.~Futer, E.~Kalfagianni, X.-S.~Lin, 
   N.~Stoltzfus, {\it The Jones polynomial and graphs on surfaces}, 
   Journal of Combinatorial Theory, Ser.B {\bf 98} (2008) 384--399;
   preprint \verb#math.GT/0605571#. 

\bibitem[DH]{DH} Y.~Diao, G.~Hetyei, {\it Relative Tutte polynomials for 
  colored graphs and virtual knot theory}, 
  Combinatorics Probability and Computing {\bf 19 }(2010) 343–-369;
  Preprint \verb#arXiv:math.CO/0909.1301#.
  
\bibitem[GT]{GT} J.~L.~Gross and T.~W.~Tucker,
   {\it Topological graph theory}, Wiley, NY, 1987.

\bibitem[K2]{Ka2} L.~H.~Kauffman, {\it A Tutte polynomial for signed graphs}, 
   Discrete Appl.~Math.~{\bf 25} (1989) 105--127.

\bibitem[K1]{Ka1} L.~H.~Kauffman, {\it New invariants in knot theory},
   Amer.~Math. Monthly~{\bf 95} (1988) 195--242.

\bibitem[K3]{Ka3} L.~H.~Kauffman, {\it Virtual knot theory}, 
   European Journal of Combinatorics {\bf 20} (1999) 663--690.

\bibitem[LZ]{LZ} S.~K.~Lando, A.~K.~Zvonkin,
   {\it Graphs on surfaces and their applications}, Springer, 2004.

\bibitem[Mof]{Mof} I.~Moffatt, {\it Knot Invariants and the 
   Bollobas-Riordan Polynomial of embedded graphs},  
   European Journal of Combinatorics {\bf 29} (2008) 95--107; 
   preprint \verb#arXiv:math.CO/0605466#.

\bibitem[MT]{MT} B.~Mohar, C.~Thomassen, {\it Graphs on Surfaces},
   The Johns Hopkins University Press, 2001.

\bibitem[Th]{Th} M.~Thistlethwaite, {\it A spanning tree expansion for
   the Jones polynomial}, Topology {\bf 26} (1987) 297--309.

\bibitem[Tr]{Tr} L.~Traldi, {\it A subset expansion of the coloured Tutte 
   polynomial}, 
   Combinatorics Probability and Computing {\bf 13 }(2004) 269–-275.

\bibitem[VT]{VT} F.~Vignes-Tourneret, {\it The multivariate signed 
   Bollobas-Riordan polynomial}, 
   Discrete Math. {\bf 309} (2009) 5968--5981; 
   preprint \verb#arXiv:math.CO/0811.1584#.

\end{thebibliography}
\end{document}